\documentclass[preprint,authoryear]{elsarticle}
\usepackage[latin1]{inputenc}
\usepackage[english]{babel}

\usepackage{amssymb}
\usepackage{amsmath,enumerate,rotate}
\usepackage{amssymb,latexsym}
\usepackage{epsfig}
\usepackage{graphicx}
\usepackage{amsthm}
\usepackage{natbib}
\usepackage{color}

\theoremstyle{plain}
\newtheorem{theorem}{Theorem}%[section]
\newtheorem*{theorem*}{Theorem}
\newtheorem{lemma}[theorem]{Lemma}

\newtheorem{proposition}[theorem]{Proposition}
\theoremstyle{remark}

 \DeclareMathOperator{\mean}{\mathbb{E}}

\newcommand\itemno[1]{(\romannumeral #1)}

\providecommand{\abs}[1]{\lvert#1\rvert}
\newcommand\Mod[1]{\vert{#1}\vert}
\renewcommand\mod[1]{\left\vert{#1}\right\vert}

\newcommand\Smallfrac[2]{\mbox{\footnotesize$\displaystyle\frac{#1}{#2}$}}
\newcommand\smallfrac[2]{\mbox{\small$\displaystyle\frac{#1}{#2}$}}

\newcommand\tsum{\textstyle\sum}
\newcommand\diff{\mathrm{d}}

\newcommand\ind{\mathbb{I}}

\newcommand\mle{{\hat \te}_n}
\newcommand\lossexpn{d_n}             % exponent of the loss function, which depends on n

\newcommand\To{\longrightarrow}

\newcommand\npost{\pi^{(n)}}
\newcommand\post{\pi_f^{(n)}}

\newcommand\BR{\mathbb{R}}

\newcommand\cA{\mathcal{A}}

\newcommand\cF{\mathcal{F}}

\newcommand\ga{\gamma}    
\newcommand\de{\delta}
  \newcommand\vep{\varepsilon}

\newcommand\la{\lambda}

\newcommand\te{\theta}

\usepackage{amssymb}
\usepackage{amsmath}
\usepackage{latexsym}
\usepackage[mathcal]{eucal}
\usepackage{mathrsfs} %per usare il carattere mathscr
\usepackage{verbatim}

\journal{Journal of Statistical Planning and Inference}

\begin{document}
\begin{frontmatter}
\title{On Bayesian Learning from Bernoulli Observations}
\author[pgb]{Pier Giovanni Bissiri\corref{cor1}} 
\ead{pier.bissiri@unimib.it}%piergiovannib@tiscali.it}%piergiovannib@yahoo.it}
\author[sgw]{Stephen G. Walker}
\ead{s.g.walker@kent.ac.uk}
\cortext[cor1]{Corresponding author}
%\maketitle
\address[pgb]{Dipartimento di Statistica, 
Università degli Studi di Milano-Bicocca,\\ 
Edificio U7, Via Bicocca degli Arcimboldi 8,\\
20126 Milano, Italy}
\address[sgw]{IMSAS, University of Kent, Canterbury, Kent, CT2 7NZ, UK\\ tel: +44(0)1227 823800}

\begin{abstract} 
We provide a reason for Bayesian updating, in the Bernoulli case, 
even when it is assumed that observations are independent and
identically distributed with a fixed but unknown parameter $\theta_0$. 
The motivation relies on the use of loss functions and asymptotics. Such a
justification is important due to the recent interest and focus on Bayesian consistency 
which indeed assumes that the observations are independent
and identically distributed rather than being conditionally independent 
with joint distribution depending on the choice of prior.

%\vspace{0.1in} \noindent {\sl Keywords:} Asymptotics; Kullback--Leibler divergence; Loss function.
\end{abstract}

\begin{keyword}
Asymptotics\sep Kullback--Leibler divergence\sep Loss function.
\end{keyword}
\end{frontmatter}

%\vspace{0.3in} \noindent {\bf 1. Introduction.} 
\section{Introduction}
The aim of this paper is to provide a straightforward and concise justification of the Bayesian
approach to updating probability beliefs, in the case of a Bernoulli sequence of random variables. 
It is then seen that the details of the result can be applied to general parametric families. 
The key is the use of a loss function combined with the notion of asymptotics. 
That is, a loss function on the space of probability distributions on $(0,1)$ is employed which uses as
information the prior knowledge and the observations. 
The general setting is made precise by appealing to obvious asymptotic requirements for
the solution to the minimization of the loss function. 
Indeed, interest in Bayesian consistency has grown in the last years. 
See, for instance, 
\cite{XiRa09} and references cited in this paper.

The use of loss functions is limitless within the world of applied sciences, 
no more so than within the decision sciences which includes 
statistics and particularly Bayesian statistics %(Berger, 1993). 
\Citep{Berger}. 
To set the scene, if $\cA$ is a set of actions, and the loss incurred is
$L(a,X)$, where $X$ is an outcome/observation or piece of information and $a\in\cA$, 
then the best choice is that $\widehat{a}$ which minimizes
$L(a,X)$. On the other hand, if there are a number of pieces of information, 
say $(X_1,\ldots,X_n)$, each of which contributes an additive loss
$L(a,X_i)$ under action $a$, then the best choice now minimizes 
the cumulative loss
$$ L(a,X_1,\ldots,X_n)=\sum_{i=1}^n L(a,X_i). $$
Such an additive style of cumulative loss would be appropriate when the $(X_i)$ 
are independent pieces of information.

We are interested in the case when 
$\cA$ is the space of probability distributions on $(0,1)$. 
This occurs if our aim is to choose a probability distribution representing beliefs about a model parameter 
$\theta$ belonging to $(0,1)$. 
In this framework, we allow $\pi(\cdot)$ 
to be the proposed representation of beliefs about $\theta$ in the case of no observations. 
The distribution $\pi$  
represents information, just as the Bernoulli observations $(X_1,\ldots,X_n)$
represent information. 
Hence, in the case $n=0$, the loss function is $l_\pi(a,\pi)$. 
Maintaining the idea of cumulative loss, we now have
\begin{equation}\label{f: loss0}
L(a,X_1,\ldots,X_n,\pi)=\sum_{i=1}^n L(a,X_i)+l_\pi(a,\pi).
\end{equation}
To this point there is little justification required; 
we are merely writing down a general loss function in order to determine a probability
distribution on $(0,1)$, where the only assumption is that the losses are additive or cumulative. 
This seems relevant when the pieces of
information are independent; 
that is, no one piece of information provides information about any of the others. 
To better indicate that $\cA$ is
a set of probability distributions, we will now replace $a$ with $\nu$. 

The first and straightforward loss function to discuss is $l_\pi(\nu,\pi)$. 
So, $l_\pi(\nu,\pi)$ is the loss when $\nu$ is the probability
measure correctly representing beliefs 
(and indeed with consistency it will end up providing correct beliefs) and $\pi$ is the proposed
probability measure representing beliefs at the outset. 
Therefore,  $l_\pi(\nu,\pi)$ can be interpreted as a loss in information. 
It is reasonable to require $\nu$ to be absolutely continuous
with respect to $\pi$. Indeed, the updated probability should be zero on every event whose
prior $\pi$ probability is zero. 
Therefore, $l_\pi(\nu,\pi)$ can be taken to be 
the $g$--divergence, introduced by 
\cite{AliSilvey} and \cite{Csiszar}, 
i.e. 
$l_\pi(\nu,\pi)=\int g(\diff\nu/\diff\pi)\,\diff\nu$, where $g$ is a convex function such that $g(1)=0$. 
Such a family of divergences is known to be a generalization of the Kullback--Leibler divergence, 
introduced by \cite{KullbackLeibler}, which is obtained taking $g(x)=-\log(x)$. 
\cite*{BissiriWalker} 
establish that among the $g$--divergences, the Kullback--Leibler divergence is 
the only one which preserves a necessary coherence property whereby the solution at stage $n$ 
serves as the prior for subsequent observations. 
Hence,  $l_\pi(\nu,\pi)$, the loss in information in using $\pi$ rather than $\nu$ is taken to be 
the Kullback--Leibler divergence. 
%$l_\pi(\nu,\pi)=\int \ln(\nu/\pi)\nu$.   

Our aim now is to ascertain how $\pi$ changes to $\nu$, 
in the light of the information $(X_1,\ldots,X_n)$, for an apparent arbitrary loss
function $L(\nu,X)$. Our form for this seems obvious in the sense that 
if we select $L(\theta,X)$ then we can merely take $L(\nu,X)=\int
l(\theta,X)\,\nu(\diff \theta)$, since the $\nu$ represents beliefs in $\theta$ 
and so $L(\nu,X)$ is understood as expected loss. Surprisingly
now, an obvious asymptotic requirement will pin down $l(\theta,X)$ precisely. 
The following can also be seen as providing an explicit answer to
a suggestion in \citet[Section 6]{Walker}  
about a possible justification of the Bayesian paradigm 
through the loss  \eqref{f: loss0} and asymptotic requirements. 
So, while loss functions are typically regarded as a subjective choice, an objective choice based on
asymptotic properties for the $\nu$ provides justification for the Bayesian learning process. 
All the proofs are deferred to the Appendix.

\section{Preliminaries}
Denote by $(X_n)_{n\geq 1}$ the sequence of observations which are i.i.d. Bernoulli with
parameter $\theta_0$. Assume they are $0\!-\!1$ random variables on a probability space 
$(\Omega, \mathscr{F}, P_\te)$ and $P_\te(X_n=1)=\te$
for each $n\geq 1$. Denote by $\pi$ the prior distribution for $\te$. 
So, $\pi$ is a probability measure on the Borel subsets of $(0,1)$. 
In the rest of the paper, it will be assumed that $\pi$ 
is absolutely continuous with respect to the Lebesgue measure $\la$ and that 
there is a version of $\diff \pi/\diff \la$ that is a continuous function.

The observations $X_1, X_2, \dotsc$ are usually considered conditionally independent 
and identically distributed given $\te$. See for example 
%Bernardo and Smith (1994). 
\cite*{Bernardo}. 
This makes possible to update the prior $\pi$ by Bayes' theorem, obtaining the posterior distribution for $\te$
\begin{equation*}%\label{f: posterior}
\npost(A):= \pi(A\mid X_1,\dotsc,X_n) =
\frac{\int_A \te^{n \mle}\, (1-\te)^{n(1-\mle)}\, \pi(\diff \te)} {
\int_{(0,1)} \te^{n \mle}\, (1-\te)^{n(1-\mle)}\, \pi(\diff \te)},\end{equation*} 
where $A$ is a Borel subset of $(0,1)$ and $\mle:=\Smallfrac{1}{n}\tsum_{i=1}^n X_i$.
 
Applying Bayes' theorem to obtain the posterior distribution, the observations are not considered independent, but conditionally independent given $\theta$. 
Since we are assuming that the observations are independent, 
we are uncomfortable with the notion of a Bayesian model which artificially 
creates a dependence between the observations. 
However, as it will be soon clear, the posterior distribution also arises as the solution of a minimization problem, 
which does not require such an assumption of %conditional independence for the observations.
dependence for the observations.

Following Section 1, we consider the loss \eqref{f: loss0} taking 
$l(\theta,X)=-\ln (P_\te(X_1=X)/P_{\theta_0}(X_1=X))$, the self--information loss function. 
So, the loss function \eqref{f: loss0} 
becomes: 
\begin{equation}\label{f: loss1}\begin{split}
L(\nu)\,:&=\, L(\nu,x_1,\dotsc,x_n,\pi)\\&= \,
-\sum_{i=1}^n \int_{(0,1)} \ln(P_\te(X_1=x_i)/P_{\theta_0}(X_1=x_i))\ \nu(\diff \te)\ +\, \ D(\nu, \pi)
%int_{(0,1)}\ \ln\bigg(\frac{\nu(\diff \te)}{\pi(\diff \te)}\bigg)\; \nu(\diff \te)
\end{split}\end{equation}
where $(x_1,\dotsc, x_n)$ is a sample drawn from $(X_1,\dotsc, X_n)$, $\nu$ 
is a probability measure on $(0,1)$ absolutely continuous
with respect to $\pi$, and $D$ denotes the  Kullback--Leibler divergence (relative entropy), i.e.
\[D(Q_1,Q_2) = \int_S \ln\bigg(\frac{\diff Q_1}{\diff Q_2}\bigg)\ \diff Q_1,\]
where $S$ is the support of $Q_1$, 
for any couple $(Q_1,Q_2)$ of probability measures such that $Q_1\ll Q_2$. 
%This is based on the loss in Section 1 with $l(\theta,X)=-\ln(P_\te(X_1=X)/P_{\theta_0}(X_1=X))$, the self--information loss function.
%The usual convention, based on continuity arguments, will be applied, taking $\ln(0)\:=\:0$. 
%$0\log(0/q)\: =\:0$ for all real $q$ and $p\,\log(p/0)\: =\: 1$ for all real non-zero $p$. 
%Hence the relative entropy assumes values in $[0,\infty]$.

%In fact, %the loss functioni in \eqref{f: loss} can be written in the following form:
Notice that the first addendum in \eqref{f: loss1} depends on the sample 
and attains its minimum when $\nu=\delta_{\mle}$, 
i.e. $\nu$ is degenerate at $\mle$, 
while the second term takes into account only the prior belief 
about $\te$ expressed by $\pi$.
It is clear that the posterior $\npost$ minimizes the loss $L$, since
\[L(\nu) =
D(\nu, \npost) - \ln\bigg(\int_{(0,1)}\, \prod_{i=1}^n P_\te(X_1=x_i)\,  \pi(\diff \te)\bigg)+\sum_{i=1}^n \ln(P_{\theta_0}(X_1=x_i)).\] We
could stop here since we have a justifiable loss function 
the solution of which is the Bayesian posterior distribution. 
However, we can work
with a more general loss function and establish through asymptotic arguments that 
the $-\log$ loss is the best in some sense.

So, an obvious and general alternative for the loss function in \eqref{f: loss1}
is %can be obtained by replacing the function $-log$ with a more general function $f$
\begin{equation}%\label{f: loss1}
\begin{split}
L_f(\nu)
:\!&= \sum_{i=1}^n \int_{(0,1)} \,f(P_\te(X_1=x_i))\ \nu(\diff \te)\\
&+\ \int_{(0,1)} \ln\bigg(\frac{\nu(\diff \te)}{\pi(\diff \te)}\bigg)\;
\nu(\diff \te),\quad \nu \ll \pi,
\end{split}
\end{equation}
where the function $-\ln(\cdot)$ has been replaced by a function $f$ from $(0,1)$ 
into the non--negative real line including $+\infty$. Clearly this is appropriate as
$f(P_\te(X_1=x))$ is either $f(\theta)$ or $f(1-\theta)$, depending on whether $x$ is 1 or 0. 
Denote by $\post$ the probability measure that is
absolutely continuous with respect to $\pi$ with density
\[\frac{e^{-n \{\mle f(\te)+(1-\mle)f(1-\te)\}}}{\int_{(0,1)} e^{-n \{\mle f(t)+(1-\mle)f(1-t)\}\, }\diff \pi(t)}.\]
%\[\frac{e^{-n \, h(\te,\mle)}}{\int_{(0,1)} e^{-n \, h(t,\mle)}\, \diff \pi(t)},\]
%where $h(x,y)\::=\:y\,f(x)\;+\;(1-y)\,f(1-x)$ for every $(x,y)$ in $(0,1)^2$.
The probability measure $\post$ minimizes $L_f$ since
%\[L_f(\nu) = D(\nu, \post) - \ln\bigg(\int_{(0,1)}\,e^{-n \{\mle f(\te)+(1-\mle)f(1-\te)\}}\,  \pi(\diff \te)\bigg)-\sum_{i=1}^n f(P_{\theta_0}(X_1=x_i)).\]
\[\begin{split}
L_f(\nu) &= D(\nu, \post) - \ln\bigg(\int_{(0,1)}\,
e^{-n \{\mle f(\te)+(1-\mle)f(1-\te)\}}\,  \pi(\diff \te)\bigg).
\end{split}\]
A referee has pointed out connections with Lagrange functions, $\post$ and the unique minimization of $L_f$.

\section{Theory}
Our aim is to properly choose the function $f$ within the class $C^1(0,1)$, apart from an additive
constant. In fact, for any real constant $c$, $\pi^{(n)}_f(\cdot)= \pi^{(n)}_{f+c}(\cdot)$.

It will be shown which conditions on $f$ are necessary and sufficient 
for the (strong) consistency of $\post$. Next, a criterion will be defined
that makes $f(\cdot)=-\ln(\cdot)$ the best choice and therefore the Bayesian posterior $\npost$ the best one.

\subsection{Consistency}
Assume that
\begin{equation}\label{f: support}\pi(\te_0-\vep, \te_0+\vep)> 0,\end{equation}
for every $\vep >0$. Since $\te_0$ is unknown, this means that only priors 
whose support is the unit interval will be considered. It will be
convenient to express the posterior in the following form:
\begin{equation}\label{f: post}
\post(A) = \frac{\int_A\ \ e^{-n \, d(t,\mle)}\, \diff \pi(t)}
{\int_{(0,1)} e^{-n \, d(t,\mle)}\, \diff \pi(t)},
\end{equation}
where
\begin{equation}\label{f: dist}
d(x,y):=%h(x,y)-h(y,y)=
\:y\,(f(x)-f(y))\;+\;(1-y)\,(f(1-x)-f(1-y))
\end{equation}
for every $(x,y)$ in $(0,1)^2$.

%Notice that if $f(\cdot)=-\ln(\cdot)$, then $d(\te_1,\te_2)=D(\mu_{\te_1},\mu_{\te_2})$, where
%\linebreak $\mu_{\te_i}=\te_i\, \delta_{\omega_1} + (1-\te_i)\,\delta_{\om_2}$ ($i=1,2$) for any two distinct elements $\om_1, \om_2$ of
%$\Om$.

%The posterior $\post$ is said to be (strongly) \emph{consistent} if for eve\[\]

\begin{proposition}\label{prop: consist1}
%Let the sequence $(X_n)_{n\geq 1}$ be Bernoulli with parameter $\te_0$ and assume
%Assume that
%\begin{equation}\label{f: support}
%\pi(\te_0-\vep, \te_0+\vep)> 0,\end{equation} for every $\vep >0$.
Let $f$ be a function of class $C^1(0,1)$ and let $d$ be a function
defined on $(0,\, 1)^2$ by \eqref{f: dist}.

Then the following facts are equivalent:
\begin{enumerate}
\item[\itemno1]
For every $\te_0$, $\vep>0$, and every prior $\pi$ satisfying \eqref{f: support},
\begin{equation}\label{f: consistency}
\post(\te_0-\vep,\te_0+\vep)\To 1,
\quad n \to \infty,\qquad P_{\te_0}-\, \text{a.s.}\end{equation}
\item[\itemno2]
For every $(\te_1,\te_2)$ in $(0,1)^2$,
\begin{enumerate}
\item[a)] $d(\te_1,\te_2)\geq 0$
\item[b)] if $\te_1 \neq \te_2$ then
$d(\te_1,\te_2)>0$.
\end{enumerate}
\item[\itemno3]
For every $0<x<1$,
\begin{enumerate}
\item[a)] $x\, f'(x)\;=\;(1-x)\,f'(1-x)$,
\item[b)] $f'(x)\,<\,0$.
\end{enumerate}
\end{enumerate}

\end{proposition}

%%%%%%%%%%%%%
In the rest of the paper, it will be assumed that $\pi$ is absolutely continuous with respect to the
Lebesgue measure and that its density is continuous on $(0,\:1)$.
The following proposition determines the rate of convergence of $\post$. 
In its statement and in the rest of the paper, 
the minimum between two real numbers $x$ and $y$ will be denoted by $x\wedge y$. 
%In the rest of the paper, the minimum and the maximum between two real numbers $x$ and $y$ will be denoted by $x\wedge y$ and $x\vee y$, respectively.
\begin{proposition}\label{prop: conv.rate}
If the hypotheses of Proposition \ref{prop: consist1} are satisfied,
%$\pi$ is absolutely continuous w.r.t. the Lebesgue measure $\la$ and $\diff \pi/\diff \la$ is continuous,
then, as $n\to \infty$,
\begin{equation}
\post((\te_0-\vep,\te_0+\vep)^c)
\asymp  \frac{1}{\sqrt{n}}\,e^{-n(\de(\vep))}
\qquad P_{\te_0}-\,\text{a.s.},
\end{equation}
where $\de(\vep):=d(\te_0-\vep, \te_0)\wedge d(\te_0+\vep, \te_0)$.

\end{proposition}

\subsection{The choice of the loss function.}
Here we study the large sample property of the posterior, and this can be
done by considering the posterior variance given $\mle$.

\begin{proposition}\label{prop: variance}
If $V_f(\mle)$ denotes the variance with respect to the distribution $\post$ and
$f$ satisfies the conditions of Proposition \ref{prop: consist1}, then
\begin{equation}\label{f: variance}
\lim_{n\to \infty}nV_f(\mle)\,=\,
-\frac{\te_0(1-\te_0)}{\te_0\,f'(\te_0)},\qquad P_{\te_0}-\text{a.s.}
\end{equation}
\end{proposition}

So, the limit depends obviously on $\theta_0$ and it is clear how:  
the Fisher information is $I(\theta_0)^{-1}\propto\theta_0(1-\theta_0)$ 
and so for larger values of $I(\theta_0)$ 
we will have a faster rate of convergence since there is more
information in the data for such $\theta_0$. 
Moreover, the information tells us how the convergence depends on $\theta_0$. 
This is amplified by the speed at which the $\mle$ converges to $\theta_0$,  
and is proportional to $\theta_0(1-\theta_0)$. 
So,  the $\theta_0(1-\theta_0)$ term in the
limit of (\ref{f: variance}) is taking account of the value of $\theta_0$. 
The other term should therefore not depend on $\theta_0$.

The reason for this is quite simple: if the $\theta_0\,f'(\te_0)$ does depend on $\theta_0$ 
then we should be able to modify $f$ so that all
$\theta_0$ obtain the largest value of $|\theta f'(\theta)|$. 
Hence, the optimal $f$ must indeed make this a constant and so we must take
$-xf'(x)=M$ for some constant $M>0$. Hence we have $f(x)=-M\ln x$.

Hence, we now just need to ascertain the reason why we should make $M=1$; 
since we have established that we must have $f(x)=-M\ln x$ and the Bayesian
learning rule is obtained precisely with $M=1$. 
Suppose the choice $\theta=\theta_1$ is chosen stubbornly so that
$\pi(\theta)=\delta_{\theta_1}(\theta)$. Hence, since $\delta_{\theta_1}$ will always represent beliefs, 
according to definitions in Section 2,
$$L(\delta_{\theta_1})=M\sum_{i=1}^n \ln\left\{\frac{P_{\theta_0}(X_1=x_i)}{P_{\te_1}(X_1=x_i)}\right\}$$
and so our expected loss for $n$ observations is
$$\bar{L}(\delta_{\theta_1})= n M D(P_{\theta_0},P_{\theta_1}).$$
We can understand that our loss up to a sample of size $n$ when fixing $\theta_1$; 
it is predicting with the wrong measure, i.e. $P_{\theta_1}$
instead of $P_{\theta_0}$ on $n$ occasions. 
So our loss is $n D(P_{\theta_0},P_{\theta_1})$, being consistent with using $D(\nu,\pi)$ in Section
1, and hence we must fix $M=1$.

%\vspace{0.2in} \noindent {\bf 4. Discussion.} 
\section{Discussion}
We have constructed a loss function for selecting an updated belief probability measure on $(0,1)$
in the light of i.i.d. Bernoulli random variables. Having started out with a general form, 
the precise function can be pinned down by appealing
to some necessary asymptotic properties. 
The consequence is that the Bayesian learning machine, in the Bernoulli case at least, can be understood
via notions of loss functions and asymptotics and while retaining the correct notion of an i.i.d. sample.

We believe that the ideas in this paper can be extended to the more general case; 
in the first instance for parametric models $f(x;\theta)$ and
subsequently for nonparametric models $f(x)$, $f\in\cF$, 
where now the decision space consists of probability measures on $\cF$.

%\subsection*{Appendix}
%\appendix
\section*{Appendix}

In order to prove Proposition \ref{prop: consist1}, the following lemma will be useful.
\begin{lemma}\label{lemma: deriv}
%Define $d$ by \eqref{f: dist}
Let $f$ be a function of class $C^1(0,\,1)$ such that
\begin{enumerate}
\item[a)] $x\, f'(x)\;=\;(1-x)\,f'(1-x)$,
\item[b)] $f'(x)\,<\,0$,
\end{enumerate}
for every $0<x<1$.
Moreover, fix $0<\te<1$ and define
\begin{equation}\label{f: dist.}
\varphi(x)\,=\,\te\,f(x)\,+\,(1-\te)\,f(1-x), \end{equation}
for every $0<x<1$.

Hence, $\varphi$ is a function of class $C^1(0,\: 1)$ such that
 \begin{equation}\label{f: loss.deriv.}
\varphi'(x)\;=\; \frac{\te - x}{1-x}\,f'(x). \end{equation}
Moreover, $\varphi$ has the second derivative at $\te$, which is equal to
\begin{equation}
\varphi''(\te)\;=\; -\frac{f'(\te)}{1-\te}. \label{f: second.deriv..}
\end{equation}

\end{lemma}

\begin{proof}
By \eqref{f: dist.},
\begin{equation}\label{f: loss.deriv1.}
\varphi'(x)\;=\;\te f'(x)\:-\:(1-\te)\,f'(1-x). \end{equation}
A combination of \eqref{f: loss.deriv1.} with (a) yields \eqref{f: loss.deriv.}.
Since $f'$ is a continuous function, \eqref{f: loss.deriv.} entails
\[\lim_{x\to\te} \frac{\varphi'(x)-\varphi'(\te)}{x-\te}\,=\,
-\frac{f'(\te)}{1-\te},\]
and \eqref{f: second.deriv..} is proved.
\end{proof}

\begin{proof}[Proof of Proposition \ref{prop: consist1}]

Let $A^c:=(0,1)\setminus A$ denote the complement of subset $A$ of $(0,1)$.

To begin with, notice that by \eqref{f: post}
\[
\post(\te_0-\vep,\te_0+\vep) =
\Bigg( 1+
\frac{\int_{(\te_0-\vep,\te_0+\vep)^c}
                                        e^{-n\, d(\te,\mle)} \pi(\diff \te)     }
                                                {\int_{(\te_0-\vep,\te_0+\vep)}
                                                e^{-n\, d(\te,\mle)} \pi(\diff \te)     }
\Bigg)^{-1}, \]
and therefore \eqref{f: consistency} is tantamount to
\begin{equation}\label{f: consist}
\lim_{n\to \infty}          \frac{\int_{(\te_0-\vep,\te_0+\vep)^c}
                                        e^{-n\, d(\te,\mle)} \pi(\diff \te)     }
                                                {\int_{(\te_0-\vep,\te_0+\vep)}
                                                e^{-n\, d(\te,\mle)} \pi(\diff \te)     }
=0  \qquad  P_{\te_0}-\text{a.s.}                   \end{equation}

Let us prove that (ii) is necessary for (i). 
To this aim, assume that (i) is true and fix $0<\te_0<1 $ and a probability measure $\pi$
satisfying \eqref{f: support}. Hence, by virtue of (i), $\pi$ must satisfy \eqref{f: consistency} as well.

Since $f$ is continuous, the function $d(\cdot,\te_1)$ is continuous as well for every $\te_1$ in $(0,1)$. 
In particular, for every
$\te_1$, it is continuous at $\te_1$ where its value is zero; i.e.
for every $k > 0$ there exists some $\vep>0$ such that
\( \inf_{\te:\,\abs{\te-\te_1}\,<\,2\vep}\,
d(\te,\te_1) > -2k \).
%for every $\te_1$.
Moreover, by the strong law of large numbers, $\abs{\mle-\te_0}<\vep$
for sufficiently large $n$ $P_{\te_0}$-a.s., and therefore
\[ \inf_{\te:\,\abs{\te-\te_0}\,<\,\vep} d(\te,\mle)\, \geq\,
\inf_{\te:\,\abs{\te-\mle}\,<\,2\vep} d(\te,\mle) > -2k
,\]
for sufficiently large $n$, $P_{\te_0}$-a.s., for every $k>0$
and some $\vep>0$.
Hence, by the dominated convergence theorem,
\begin{equation}\label{f: piribi}
\begin{split}
\lim_{n\to\infty} e^{-2nk}\;
&\int_{(\te_0-\vep,\te_0+\vep)}
e^{-n\, d(\te,\mle)} \pi(\diff \te)\\   &=
\int_{(\te_0-\vep,\te_0+\vep)} \lim_{n\to\infty}
e^{-n\, (d(\te,\mle)+2k)} \pi(\diff \te)
=0
\qquad P_{\te_0}-\text{a.s.},
\end{split}\end{equation}
for every $k>0$ and for some $\vep>0$.
%Therefore, for some constant $c$,
%\begin{equation}\label{f: piribi}\int_{(\te_0-\vep,\te_0+\vep)}
%e^{-n\, d(\te,\mle)} \pi(\diff \te)\leq c\,e^{nk},\ \end{equation}eventually.

%Since $d$ is continuous, it is bounded on any compact subset $K$ of $[0,1]$.
Now assume that there is some $M$ such that $d(\te,\te_0)\leq M$ for every $\te$ belonging to some set $C$ 
such that $\pi(C\setminus \{\te_0\})>0$.
Take $\vep\geq 0$
small enough so that $\pi((\te_0-\vep,\te_0+\vep)^c \cap C) >0$.
Notice that, by the strong law of large numbers and by continuity of $d(\te, \cdot),$ for every $\te\in(0,1)$,
$d(\te,\mle)-d(\te,\te_0)\,<\,M$ for sufficiently large $n$,
$P_{\te_0}$-a.s.
 Hence, for every $\te\in C$, $d(\te,\mle)\,<\,2M$ for sufficiently large $n$,
$P_{\te_0}$-a.s.
and,
 by Fatou's lemma,
\begin{equation}\label{f: piribo}\begin{split}
\liminf_{n\to\infty} e^{2nM}\;
&\int_{(\te_0-\vep,\te_0+\vep)^c}
e^{-n\, d(\te,\mle)} \pi(\diff \te)\\
   &\geq
\int_{(\te_0-\vep,\te_0+\vep)^c\cap C} \liminf_{n\to\infty} e^{n\, (2M-d(\te,\mle))}\, \pi(\diff \te) =\infty,
\qquad P_{\te_0}-\text{a.s.}
\end{split}\end{equation}
%Therefore, for some $\vep>0$,
%\begin{equation}\label{f: piribo}
%\int_{(\te_0-\vep,\te_0+\vep)^c}
%e^{-n\, d(\te,\mle)} \pi(\diff \te)\geq \,e^{-nM},
%\end{equation}
%for all $n$ but finitely many.

So, combining \eqref{f: piribi} and \eqref{f: piribo}, one notices that
\begin{equation}
\frac{\int_{(\te_0-\vep,\te_0+\vep)^c} e^{-n\, d(\te,\mle)} \pi(\diff \te)} {\int_{(\te_0-\vep,\te_0+\vep)} e^{-n\, d(\te,\mle)} \pi(\diff \te)}
\geq e^{-2n(M+k)}
\end{equation}
holds for sufficiently large $n$, $P_{\te_0}$-a.s.,
for every $k>0$ and some $\vep>0$. Since  \eqref{f: consist} holds true for every $\vep\geq 0$, then $M \geq -k$ for any real positive
number $k$, i.e. $M \geq 0$. So, $M\geq 0$ whenever $d(\te,\te_0)\leq M$ with  positive $\pi$-probability.
%$\te\neq \te_0$, with positive $\pi$-probability.
Therefore, $d(\te,\te_0)\geq 0,$ $\pi$-a.s.
Hence,
\begin{equation}\label{f: piribi2}
\begin{split}
\lim_{n\to\infty}
&\int_{(\te_0-\vep,\te_0+\vep)}
e^{-n\, d(\te,\mle)} \pi(\diff \te)\\   &=
\int_{\{\te\in(\te_0-\vep,\te_0+\vep):\,d(\te,\te_0)\,>\,0\}}
 \lim_{n\to\infty}
e^{-n\, d(\te,\mle)} \pi(\diff \te)\\
&\phantom{\lim}+\,
\lim_{n\to\infty}
\int_{\{\te\in(\te_0-\vep,\te_0+\vep):\,d(\te,\te_0)\,=\,0\}}
e^{-n\, d(\te,\mle)} \pi(\diff \te)\\
&=\pi\{\te\in(\te_0-\vep,\te_0+\vep):\,d(\te,\te_0)=0\}\\
&\leq \pi(\te_0-\vep,\te_0+\vep)\}
\qquad\quad P_{\te_0}-\text{a.s.},
\end{split}\end{equation}
holds true $P_{\te_0}$-a.s., by dominated convergence theorem. Assume there is some $M$ such that $d(\te,\te_0)\leq M$,
for every $\te$ belonging to some set $D$ with  positive $\pi$-probability and such that $\te_0\notin D$.
Hence, combining \eqref{f: piribi2} and \eqref{f: piribo}, one notices that
\begin{equation}
\frac{\int_{(\te_0-\vep,\te_0+\vep)^c} e^{-n\, d(\te,\mle)} \pi(\diff \te)} {\int_{(\te_0-\vep,\te_0+\vep)} e^{-n\, d(\te,\mle)} \pi(\diff \te)}
\geq \pi(\te_0-\vep,\te_0+\vep)\, e^{-2nM}
\end{equation}
holds true for sufficiently large $n$, $P_{\te_0}$-a.s., for sufficiently small $\vep>0$. Since  \eqref{f: consist} holds true for every
$\vep\geq 0$ together with \eqref{f: support}, then $M$ must be positive.
So, $M > 0$ whenever $d(\te,\te_0)\leq M$ for every $\te \neq \te_0$
with positive $\pi$-probability.
%$\te\neq \te_0$, with positive $\pi$-probability.
Therefore, $d(\te,\te_0)>0,$ for every $\te \neq \te_0$ $\pi$-a.s.
Since this is true for every $\te_0$ and every $\pi$ whose support is the unit interval, (ii.b) must hold.
Notice that (ii.b) trivially entails (ii.a) since $d(\te,\te)=0$ for every $\te$ by definition of $d$.

At this point, it will be proved that (ii) implies (iii).
%and (iii) are equivalent.
To this aim, define $\varphi(x) := d(x,\te)$ for a fixed $0<\te<1$.
%and notice that \begin{equation}\label{f: loss.deriv1}\varphi'(x)\;=\;\te\,f'(x)\:-\:(1-\te)\,f'(1-x). \end{equation}
Since $d(\te,\te)=0$, condition (ii) is tantamount to say that the function
$\varphi$ has an absolute minimum at $x=\te$ for any $\te$.
Therefore,
if (ii) is in force,
$\varphi'(\te) = 0$ must be true for every $\te$ in the unit interval and
condition (iii.a) follows.

By Lemma \ref{lemma: deriv}, (iii.a) entails
 \begin{equation}\label{f: loss.deriv}
\varphi'(x)\;=\; \frac{\te - x}{1-x}\,f'(x), \end{equation}
where $\varphi'$ is a continuous function since $f'$ is so. Since $\te$ is an
absolute minimum point for $\varphi$, there is $\delta >0$ 
such that $\varphi'(x) > 0$ if $\te < x \leq \te+\delta$ and $\varphi'(x) < 0$  if
$\te -\delta \leq x < \te$. 
This is tantamount to say that $f'(x)<0$ for every $x$ in $(\te-\de,\te)\cup (\te,\te+\de)$ for some $\de$. Since
this must hold for every $\te$, condition (iii.b) follows.

%On the other hand, it is clear if (iii) is true, then \eqref{f: loss.deriv} holds and therefore (ii) is also true.

Finally, it will be shown that (iii) is sufficient for (i).
To this aim, notice that if (iii) holds then \eqref{f: loss.deriv} is also in force
by Lemma \ref{lemma: deriv}
and therefore $d(\cdot,\te)$ is (strictly) decreasing on $(0,\te)$
and (strictly) increasing on $(\te,1)$, for every $\te$.
Therefore, for every
$\vep>0$ and every $0<\te_1<1$, $d(\te, \te_1) > \de(\vep, \te_1)/2$ 
if $\mod{\te-\te_1}>\vep$ and $\de(\vep,\te_1)$ denotes
$d(\te_1-\vep, \te_1)\wedge d(\te_1+\vep, \te_1)$.
%In what follows, $x\wedge y$ and $x\vee y$ will stand for the minimum and the maximum between $x$ and $y$, respectively.
Applying dominated convergence theorem, this entails that
\begin{equation}\label{f: proof.numerator}%\label{f: proof.suff1}
\begin{split}
\lim_{n\to\infty}\, &e^{n\de(\vep,\mle)/2}\;
\int_{(\te_0-\vep,\te_0+\vep)^c}
e^{-n\, d(\te,\mle)}\; \pi(\diff \te)\\  &=
\int_{(\te_0-\vep,\te_0+\vep)^c}
\lim_{n\to\infty}\,
e^{-n\,(d(\te,\mle)-\de(\vep,\mle)/2)}\; \pi(\diff \te)
=0.
%\\ &=\pi((\te_0-\vep,\te_0+\vep)^c).
\end{split}
\end{equation}

%Hence,
%\begin{equation}\label{f: proof.numerator}
%\int_{(\te_0-\vep,\te_0+\vep)^c}
%e^{-n\, d(\te,\mle)}\; \pi(\diff \te)
%\leq c_1\pi((\te_0-\vep,\te_0+\vep)^c)
%\, e^{-n\de(\vep)},\end{equation} holds eventually for some $c_1>0$.

By continuity of $d(\cdot,\mle)$ at $\mle$,
for every $\eta>0$ there exists $\ga$ such that
$d(\te,\mle)<\eta$ if $\abs{\te-\mle}<2\ga$ and by the strong law of large
numbers $\abs{\mle-\te_0}<\ga$ for sufficiently large $n$, $P_{\te_0}$-a.s.
Therefore $d(\te,\mle)<\eta$ if $\abs{\te-\te_0}<\ga$
for sufficiently large $n$, $P_{\te-0}$-a.s.,
and by Fatou's lemma,
\begin{equation}\label{f: proof.denominator}%\label{f: proof.suff1}
\begin{split}
\liminf_{n\to\infty}\, e^{n\eta}\;
&\int_{(0,1)}
e^{-n\, d(\te,\mle)}\; \pi(\diff \te)\\ \geq
&\int_{(\te_0-\ga,\te_0+\ga)}
\liminf_{n\to\infty}\,
e^{n\,(\eta- d(\te,\mle))}\; \pi(\diff \te) = \infty,
\end{split}
\end{equation}
for every $\eta >0$ and some $\ga>0$, $P_{\te-0}$-a.s..
%Therefore,\begin{equation}\label{f: proof.denominator}\int_{(0,1)}
%e^{-n\, d(\te,\mle)}\; \pi(\diff \te) \geq c_2 e^{-n\eta},\end{equation}
%holds eventually for every $\eta>0$ and some $c_2>0$.

Combining \eqref{f: proof.numerator} and \eqref{f: proof.denominator}, one obtains that
\[
\lim_{n\to\infty}\; e^{n(\de(\vep,\mle)/2-\eta)} \:
\frac{\int_{(\te_0-\vep,\te_0+\vep)^c}
e^{-n\, d(\te,\mle)}\; \pi(\diff \te)}
{\int_{(0,1)}
e^{-n\, d(\te,\mle)}\; \pi(\diff \te)}=0,
\]
for every $\eta,\vep>0$.
Taking $\eta<\delta(\vep,\te_0)/2$, (i) follows.

\end{proof}

By the strong law of large numbers, there exists a Borelian subset 
$B$ of $\{0,1\}^\infty$ with $P_{\te_0}$-probability one such that 
$\mle(x_1,\dotsc,x_n)=\Smallfrac{1}{n}\tsum_{i=1}^n x_i$ converges to $\te_0$ for all sequences $(x_n)_{n\geq 1}$ 
belonging to $B$.
In the rest of this appendix, $\mle$ will stand 
for $\mle(x_1,\dotsc,x_n)$ and we shall always assume that $(x_n)_{n\geq 1}$ belongs to $B$.

In order to prove Proposition \ref{prop: conv.rate}
and Proposition \ref{prop: variance},
the following lemmas are useful.

\begin{lemma}\label{lemma: interv.}
If $(g_n)_{n\geq 0}$ is a sequence of non-negative functions on $(0,\: 1)$
dominated by an integrable function,
then for every $\delta>0$ there are $\eta_1,\, c_0>0$ such that
\begin{equation}
\int_{(\mle-\delta,\,\mle+\delta)^c}\,
e^{-nd(t, \mle)}\, g_n(t)\,\diff t \leq c_0\,e^{-n\eta_1}
\label{f: int.dist.}
\end{equation}
for sufficiently large $n$.

Moreover, if $(c_n)_{n\geq 1}$ is a sequence converging to a positive number,
$(g^*_n)_{n\geq 1}$ is a sequence of integrable functions on $\BR$, and
\begin{equation}\label{f: integral...}
\int_{(-\infty,\,\infty)}\,
e^{-c_n(t-\mle)^2/2}\,g^*_n(t)\, \diff s \,<\,c,
\end{equation}
for some real constant $c$ and for sufficiently large $n$, then
\begin{equation}
\int_{(\mle-\delta,\,\mle+\delta)^c}\,
e^{-n\,c_n(t-\mle)^2/2}\,g^*_n(t-\mle)\, \diff t \leq k\,e^{-n\eta_2},
\label{f: int.gauss.}
\end{equation}
for some $k,\,\eta_2>0$ and for sufficiently large $n$.
\end{lemma}
In the rest of the paper, the maximum between two real numbers 
$x$ and $y$ will be denoted by $x\vee y$. 
\begin{proof}
Let $\chi_n$
be a nonnegative and differentiable function on $(a,\,b)$
($-\infty\leq a,b\leq \infty$)
with an unique absolute minimum at $\mle$
and such that %$\chi_n'$ is positive on $(0,\:b)$ and negative on $(a,\:0)$.
\begin{equation}\label{f: condiz.funz}\begin{split}
\chi'_n(t)&<0 \qquad \text{if}\quad a<t< \mle \\
\chi'_n(t)&>0 \qquad \text{if}\quad \mle <t<b.
\end{split}\end{equation}
Hence, if $\delta>0$
and 
$\Mod{t-\mle}>\delta$ then $\chi_n(t)>\eta_n(\delta)$, where $\eta_n(\delta):=\chi_n(\mle-\delta)\vee\chi_n(\mle+\delta)$.
Notice that $-n\,\chi_n(t) \,\leq\, -(n-1)\eta_n(\delta)\,-\,\chi_n(t)$
if $\Mod{t-\mle}>\delta$ and $n>1$.
Therefore, given some sequence of measures $(\mu_n)_n$ on
the Borelian subsets of $(a,\: b)$,
\begin{equation}\label{f: integral..}
\int_{\{\mod{t}\,>\,\delta\}}
e^{-n\chi_n(t)}\mu_n(\diff t) <
e^{-(n-1)\eta_n(\delta)}\int_{(a,\,b)}e^{-\,\chi_n(t)}\mu_n(\diff t). \end{equation}

Taking $\chi_n(t)=d(t, \mle),$ $a=0,$ $b=1$,
\eqref{f: condiz.funz} holds true. 
Moreover, the integral
\begin{equation}\label{f: integr.}
\int_{(a,\,b)}e^{-\,\chi_n(t)}\mu_n(\diff t)
\end{equation}
is less than a constant, if $\diff\mu_n/\diff \la =g_n$ and $\la$ is the Lebesgue measure.
In fact, $\chi_n$ is nonnegative and $g_n$ dominated.
Since
$\eta_n(\delta)=d(\mle-\delta, \mle)\vee d(\mle+\delta, \mle)$ converges to a positive constant
by the strong law of large numbers,
\eqref{f: integral..} yields \eqref{f: int.dist.}.

If $\chi_n(t)=c_n(t-\mle)^2/2,$ $a=-\infty,$ $b=\infty$, then
\eqref{f: condiz.funz} is satisfied.
Moreover, if $\diff\mu_n/\diff \la(t)=g^*_n(t-\mle)$, then the integral
\eqref{f: integr.} turns out to be equal to the integral in \eqref{f: integral...}.
Therefore, \eqref{f: int.gauss.} follows from
\eqref{f: integral..}. 
\end{proof}

\begin{lemma}\label{lemma: g_n}
Let $(g_n)_{n\geq 1}$ and $(g^*_n)_{n \geq 1}$ be two sequences of
nonnegative, continuous and integrable functions
defined on $(0,\: 1)$ and $\BR$, respectively, and such that
$g_n(t)\sim g_n^*(t)$ as $t\to \mle$, for every $n\geq 1$.

Let $\lossexpn(t)$ stand for  $d(t,\mle)$, and denote:
\begin{align}\label{f: IN.}
I_n\,:\!&=\, \int_{(0,1)}\,  e^{-n\lossexpn(t)}\,g_n(t)\, \diff t, \\
\label{f: IN}
I_n(x):\!&=
\int_{(-\infty,\,\infty)}\,
e^{-n(\lossexpn''(\mle)-x)(t-\mle)^2/2}\,g^*_n(t)\, \diff t,
\end{align}

Assume that
\begin{equation}\label{f: hyp..}
\lim_{n\to\infty} e^{-c\,n}/(I_n(x)-I_n(y))=0,\end{equation}
for every $c>0,\,$ and every $x,\,y$ belonging to some neighborhood of zero.
Moreover, let \eqref{f: int.gauss.} hold with $c_n=\lossexpn''(\mle)$,
for some positive constants $k, \eta_2$.

Therefore, $I_n\sim I_n(0)$ as $n\to \infty$.
\end{lemma}

\begin{proof}
This proof will be based on the Laplace method. See, for instance, 
%\cite{deBruijn}, pages 63--65. 
\citet[pp. 63--65]{deBruijn}. 
His results do not precisely fit our needs
and therefore we have to prove this lemma starting from scratch.

Recall that $\lossexpn(\mle)=0$. Moreover,
by hypothesis, $\lossexpn$ has a unique minimum at $\mle$ so that
$\lossexpn'(\mle)=0$. By Taylor's theorem, for each $n\geq 1$ and each $\vep>0$ there exists $\delta_n>0$ such that
if  $\mod{t-\te_0}<\delta_n$ then
\begin{equation}\label{f: taylor}
\Mod{\lossexpn(t)- \smallfrac{1}{2}\,
\lossexpn''(\mle)\,(t-\mle)^2}
\,<\,\vep\,(t-\mle)^2.
\end{equation}

It will be useful to observe that $\delta_n$ can be taken constant for sufficiently large $n$.
In order to show this fact, define
\begin{equation*}%\begin{split}
\psi_n(t)\,:=\,\lossexpn(t)\,-
\,\Smallfrac{1}{2}\lossexpn''(\mle)(t-\mle)^2,
%\\&=\,\mle f(t)\,+\,(1-\mle)f(1-t) \,+\,\Smallfrac{1}{2}(t-\mle)^2f'(\mle)/(1-\mle).
%\end{split}
\end{equation*}
so that $\psi_n(\mle)=\psi'_n(\mle)=\psi_n''(\mle)=0$.
By \eqref{f: loss.deriv.} and \eqref{f: second.deriv..},
\begin{equation}\label{f: psi}\begin{split}
\psi'_n(t)\,:\!&=\,\lossexpn'(t)\,-\,\lossexpn''(\mle)(t-\mle)\\
&=\,
(\mle -t)\,\Big(\Smallfrac{f'(t)}{1-t}-\Smallfrac{f'(\mle)}{1-\mle}\Big).
\end{split}\end{equation}

Recall that $0<\te_0<1$ and fix $0<\ga<(1-\te_0)\wedge\te_0$.
By hypothesis, the function $f'(t)/(1-t)$ is continuous over the compact set  $[\te_0-\ga,\te_0+\ga]$ and therefore is uniformly continuous over
that interval.
Moreover, recall that by the strong law of large numbers, 
$\mle$ belongs to $[\te_0-\ga,\te_0+\ga]$ if $n\geq N$ for some $N$. Hence,
by \eqref{f: psi} for every $\vep>0$
there exists $\delta>0$ such that
\begin{equation}\label{f: taylor1.}
\Mod{\psi'_n(t)/(t-\mle)}=\Mod{f'(t)/(1-t)-f'(\mle)/(1-\mle)}<\vep
\end{equation}
if $\Mod{t-\mle}\, <\,\delta$ and $n\geq N$.
By Lagrange's mean value theorem,
\begin{equation}\label{f: taylor2.}
\psi_n(t)\,=\,\psi_n(t)-\psi_n(\mle) \,=\,(t-\mle)\psi'_n(s)\end{equation}
for some $s$ between $t$ and $\mle$.
Combining \eqref{f: taylor1.} with \eqref{f: taylor2.}, one obtains
\[
\Mod{\psi_n(t)}\,=\, \Mod{(t-\mle)\,\psi'_n(s)}
\,\leq\,\vep \Mod{t-\mle}\Mod{s-\mle}.\]
Since \(\Mod{s-\mle}\leq\Mod{t-\mle}\), \eqref{f: taylor} holds true 
for every  $t\in (\te_0-\delta,\te_0+\delta)$ and every $n\geq N$.

For every $n\geq 1$, $g_n(t)\sim g_n^*(t)$ as $t\to \mle$, by hypothesis. Hence, the function
\[t\longrightarrow \ind_{\{\mle\}^c}(t)\,g_n(t)/g^*_n(t)\,+\,\ind_{\{\mle\}}(t)\]
is continuous on the compact set $[\te_0-\ga,\te_0+\ga]$ and therefore uniformly continuous on that set. 
For this reason, for each
$\vep>0$ there is $\delta>0$ such that
\begin{equation}\label{f: funz.g}
g^*_n(t)(1-\vep)\, \leq\, g_n(t) \,\leq\, g^*_n(t)(1+\vep)
\end{equation}
if $\Mod{t-\mle}<\delta$ and $n$ is sufficiently large.

At this stage, fix $\vep$ belonging to $\big(\:0,\, -f'(\te_0)/\{4(1-\te_0)\}\:\big)$ so that
\begin{equation}\label{f: epsilon}
\vep<-f'(\mle)/\{3(1-\mle)\}=\lossexpn''(\mle)/3
\end{equation}
for $n\geq M$ and some $M\geq N$. 
Moreover, take $\delta>0$ small enough so that \eqref{f: taylor} and \eqref{f: funz.g} are both satisfied.

Decompose the integral $I_n$ defined by \eqref{f: IN.} in the following way:
\[
I_n
=  \int_{(\mle-\delta,\,\mle+\delta)}\, e^{-n\lossexpn(t)}\, g_n(t)\, \diff t\,
+\,\int_{(\mle-\delta,\,\mle+\delta)^c}\,  e^{-n\lossexpn(t)}\,g_n(t)\, \diff t.
\]
The first term can be bounded by \eqref{f: taylor}, the second one by \eqref{f: int.dist.}, obtaining
\begin{equation*}\begin{split}
I_n&\geq \int_{(\mle-\delta,\,\mle+\delta)}\,
e^{-n(\lossexpn''(\mle)+2\vep)(t-\mle)^2/2}\,g_n(t)\, \diff t\\
I_n&\leq\,
\int_{(\mle-\delta,\,\mle+\delta)}\,
e^{-n(\lossexpn''(\mle)-2\vep)(t-\mle)^2/2}\,g_n(t)\, \diff t
+ c_0e^{-n\eta_1},
\end{split}\end{equation*}
which in virtue of \eqref{f: funz.g} becomes
\begin{align}
I_n&\geq (1-\vep)\int_{(\mle-\delta,\,\mle+\delta)}\,
e^{-n(\lossexpn''(\mle)+2\vep)(t-\mle)^2/2}\,g^*_n(t)\, \diff t    \label{f: In1} \\
I_n&\leq\,(1+\vep)
\int_{(\mle-\delta,\,\mle+\delta)}\,
e^{-n(\lossexpn''(\mle)-2\vep)(t-\mle)^2/2}\,g^*_n(t)\, \diff t
+ c_0e^{-n\eta_1}.   \label{f: In2}
\end{align}
By hypothesis, \eqref{f: int.gauss.} holds true
with $c_n\,=\,\lossexpn''(\mle)$,
for some positive constants $k, \eta_2$.
%for some real constant $c$.
Therefore,
\eqref{f: In1} becomes
\begin{equation}\label{f: In1.}
I_n\geq (1-\vep)\int_{(-\infty,\,\infty)}\,
e^{-n(\lossexpn''(\mle)+2\vep)(t-\mle)^2/2}\,g^*_n(t)\, \diff t
-(1-\vep)k\,e^{-n\eta_2}.
\end{equation}
Recalling \eqref{f: IN},
%where $\mod{x}< f'(\te_0)/\{4(1-\te_0)\}$,
the combination of \eqref{f: In1.} and \eqref{f: In2} yields
\begin{equation}\label{f: In}
(1-\vep)\,I_n(-2\vep)\,-\,(1-\vep)ke^{-n\eta_2}
\leq I_n\leq (1+\vep)I_n(2\vep)\, +\,c_0\, e^{-n\eta_1}
\end{equation}
for $n$ sufficiently large.
By \eqref{f: hyp..}, if $n$ is sufficiently large, then
\[\begin{split}
e^{-\eta_1 n}\,&<\,(1+\vep)(I_n(3\vep)-I_n(2\vep))/c_0\\
e^{-\eta_2 n}\,&<\,(I_n(-2\vep)-I_n(-3\vep))/k,
\end{split}\]
being $I_n(x)$ an increasing function of $x$.
Therefore, by \eqref{f: In},
\begin{equation}\label{f: In!}
(1-\vep)\,I_n(-3\vep)\,\leq\, I_n\,\leq\, (1+\vep)\,I_n(3\vep)
\end{equation}
holds true for sufficiently large $n$.
The number $\vep$ being arbitrary, it follows that
$I_n \sim I_n(0)$ as $n\to \infty$

\end{proof}

\begin{lemma}\label{lemma: var.norm}
Let $\lossexpn(t)$ stand for $d(t,\,\mle)$.
If the hypotheses and the conditions of Proposition \ref{prop: consist1} hold true,
$p$ is an integrable, nonnegative and continuous function on $(0,\: 1)$ and
$0\,<\,a\,<\,\te_0\,<\,b\,<\,1$, then
\begin{align}
&\int_{(0,\ 1)}\, p(t)\, e^{-n\,\lossexpn(t)}\, \diff t
\sim \frac{\sqrt{2\pi}\,p(\mle)}{\sqrt{n\,\lossexpn''(\mle)}},    \label{f: Den}  \\
\label{f: VAR}
&\int_{(0,1)}\, p(t)\,e^{-n \, \lossexpn(t)}\,(t-\mle)^2\, \diff t
\sim \sqrt{2\pi}\,p(\mle)\, \{n\lossexpn''(\mle)\}^{-3/2}             \\
\label{f: VARR}
&\int_{(0,1)}\, p(t)\,e^{-n \, \lossexpn(t)}\,t^2\, \diff t
\sim
\sqrt{2\pi}\,p(\mle)\, \{n\,\lossexpn''(\mle)\}^{-1/2}\;
    (\mle^2+\{n\,\lossexpn''(\mle)\}^{-1}),\\
&\int_{(b,\: 1)} p(s)\, e^{-n\,\lossexpn(s)}\, \diff s \sim
\frac{1}{n}\, \frac{p(b)}{\lossexpn'(b)}\, e^{-n\,\lossexpn(b)},\label{f: Num2} \\
&\int_{(0,\: a)} p(s)\, e^{-n\,\lossexpn(s)}\, \diff s \sim
-\frac{1}{n}\, \frac{p(a)}{\lossexpn'(a)}\, e^{-n\,\lossexpn(a)}, \label{f: Num1}
\end{align}
as $n\to\infty$.
\end{lemma}

\begin{proof}%[Proof of Lemma \ref{lemma: var.norm}]
Lemma \ref{lemma: g_n} will be applied to prove
\eqref{f: Den}, \eqref{f: VAR} and \eqref{f: VARR}.
Three cases will be considered:
\begin{enumerate}
    \item[Case A)] $g_n(t)=p(t)$,\quad $g^*_n(t)=p(\mle)$;
    \item[Case B)] $g_n(t)=p(t)(t-\mle)^2$,\quad $g^*_n(t)=p(\mle)(t-\mle)^2$;
    \item[Case C)]$g_n(t)=p(t)t^2$,\quad $g^*_n(t)=p(\mle)t^2$;
\end{enumerate}
Notice that the integral $I_n(x)$ defined by \eqref{f: IN} is finite
if $x\,<-\,f'(\te_0)/\{2(1-\te_0)\}$
and $n$ is sufficiently large.
In fact, by the strong law of large numbers, this entails that
$x\,<-\,f'(\mle)/(1-\mle)$ for sufficiently large $n$,
and therefore, by \eqref{f: second.deriv..} in Lemma \ref{lemma: deriv},
$\lossexpn''(\mle)-x$ is positive.

If $x\,<-\,f'(\te_0)/\{2(1-\te_0)\}$, then
\[I_n(x)=\sqrt{\frac{2\pi}{n(d_n''(\mle)-x)}}\,\mean(g_n^*(W_n)),\]
where $W_n$ is a Gaussian random variable with mean $\mle$ and variance $\{n(d''(\mle)-x)\}^{-1}$.
Therefore,
\begin{enumerate}
    \item[Case A)]
    %$I_n(x)\,=\,\frac{\sqrt{2\pi}\,p(\mle)}{\sqrt{n\,(\lossexpn''(\mle)-x)}}$,
    $I_n(x)\,=\sqrt{2\pi}\,p(\mle)\, \{n(\lossexpn''(\mle)-x)\}^{-1/2}$,
    \item[Case B)]
    $I_n(x)\,=\sqrt{2\pi}\,p(\mle)\, \{n(\lossexpn''(\mle)-x)\}^{-3/2}$,
        \item[Case C)]
    $I_n(x)\,=\sqrt{2\pi}\,p(\mle)\, \{n(\lossexpn''(\mle)-x)\}^{-1/2}\;
    (\mle^2+\{n(d''(\mle)-x)\}^{-1})$,
\end{enumerate}
if $x\,<-\,f'(\te_0)/\{2(1-\te_0)\}$ and $n$ is sufficiently large.

In all three cases,
\eqref{f: hyp..} holds true if $c>0,\,$ $\mod{x},\mod{y}< - f'(\te_0)/\{2(1-\te_0)\}$.
Moreover, the integral in \eqref{f: integral...} with
$c_n\,=\,n\,\lossexpn''(\mle)$ is equal to
\begin{enumerate}
    \item[Case A)] %$\frac{\sqrt{2\pi}\,p(\mle)}{\sqrt{n\,(\lossexpn''(\mle))}}$
        $\sqrt{2\pi}\,p(\mle)\, \{n\,\lossexpn''(\mle)\}^{-1/2}$,
    \item[Case B)]
                $p(\mle)\,\sqrt{2\pi}\, \{n\,\lossexpn''(\mle)\}^{-3/2}$,
    \item[Case C)]
    $\sqrt{2\pi}\,p(\mle)\, \{n\,\lossexpn''(\mle)\}^{-1/2}\;
    (\mle^2+\{n\,\lossexpn''(\mle)\}^{-1})$,
\end{enumerate}
which converge to zero by continuity of $\lossexpn''$ and $p$,
and by the strong law of large numbers.
Therefore, \eqref{f: integral...} is satisfied in all three cases if
$n$ is sufficiently large. This allows us to apply Lemma \ref{lemma: interv.} and to obtain that
\eqref{f: int.gauss.} holds true for some positive constants $k,$ $\eta_2$.

Since \eqref{f: hyp..} and \eqref{f: int.gauss.} hold for some positive constants $k, \eta_2$,
Lemma \ref{lemma: g_n} can be applied and
\eqref{f: Den}, \eqref{f: VAR} and \eqref{f: VARR} are proved.

%\[J_n:= \int_{(-\infty,\,\infty)}\,
%e^{-n(\lossexpn''(\mle))(t-\mle)^2/2}\,g^*_n(t)\, \diff t.\]

At this stage, our aim is to prove \eqref{f: Num2} and \eqref{f: Num1}.
The Laplace method will be used again.
By continuity of the function $p$,
for each $\vep>0$ there is $\delta>0$ such that
\begin{equation}\label{f: funz.g.}
p(b)(1-\vep)\, \leq\, p(t) \,\leq\, p(b)(1+\vep)
\end{equation}
if $b\leq t<b+\delta$.

Since $f$ is continuous, the functions
\begin{align*}t&\longrightarrow f(t)-f(b)-(t-b)f'(b) \\
t&\longrightarrow f(1-t)-f(1-b)-(b-t)f'(1-b)  \end{align*}
are continuous at $b$ and at $1-b$. Therefore, for a given $\vep>0$ we can fix $\delta>0$ such that
\begin{equation*}%\label{f: ABC}
\begin{split}
&\mod{\ d(t,\,s)\,-\,d(b,\,s)\,-\,(t-b)\,\frac{\partial d(x,\,s)}{\partial x}\bigg\vert_{x=b}\ }\\
&\phantom{XXXX}\leq s\,\mod{(f(t)-f(b)-(t-b)f'(b))}\\
&\phantom{XXXXXX}+\,(1-s)\mod{f(1-t)-f(1-b)-(b-t)f'(1-b)}<\vep
\end{split}\end{equation*}
holds for every $t\in (b,\,b+\delta)$ for and every $s\in (0,1)$.

Hence,
for a given $\vep>0$, we can fix $\delta>0$ such that
\begin{equation}\label{f: Taylor}
\mod{\lossexpn(t)-\lossexpn(b)-(t-b)\,\lossexpn'(b)}\leq \vep
\end{equation}
holds true together with \eqref{f: funz.g.} for every $t\in (b,\:b+\delta)$.

Denote
\[
J_n :=
\int_{(b,\:1)}e^{-n\lossexpn(t)}p(t)\diff t.
\]
Since $\lossexpn$ is increasing on $(b+\delta,\:1)\subset (\te_0,\:1)$,
\begin{equation}\begin{split}%\label{f: let's integrate!}
 J_n&=
\int_{(b,\:b+\delta)}e^{-n\lossexpn(t)}p(t)\diff t +
\int_{(b+\delta,\:1)}e^{-n\lossexpn(t)}p(t)\diff t \\
&\leq
\int_{(b,\:b+\delta)}e^{-n\lossexpn(t)}p(t)\diff t +
k\,e^{-n\lossexpn(b+\delta)}
\end{split}\end{equation}
where $k=\int_{(b+\delta,\: 1)}p(t)\diff t$. Therefore, we can write
\[
\int_{(b,\:b+\delta)}e^{-n\lossexpn(t)}p(t)\diff t \leq
J_n
\leq
\int_{(b,\:b+\delta)}e^{-n\lossexpn(t)}p(t)\diff t +
k\,e^{-n\lossexpn(b+\delta)},
\]
which yields, by \eqref{f: funz.g.} and \eqref{f: Taylor}, for sufficiently large $n$,
\begin{equation*}\begin{split}
J_n\, &\geq\,
(p(b)-\vep)\,
e^{-nd_n(b)}\,
\int_{(b,\:b+\delta)}e^{-n(\lossexpn'(b)+\vep)(t-b)}\diff t \\
J_n\, &\leq \,
(p(b)+\vep)\,e^{-nd_n(b)}\,
\int_{(b,\:b+\delta)}e^{-n(\lossexpn'(b)-\vep)(t-b)}\diff t +
k\,e^{-n\lossexpn(b+\delta)},
\end{split}\end{equation*}
that is
\begin{equation}\label{f: let's integrate !!}
(p(b)-\vep)\,\bar{J}_n(-\vep)
\leq J_n \leq
(p(b)+\vep)\, \bar{J}_n(\vep)\,+\,
k\,e^{-n\lossexpn(b+\delta)},
\end{equation}
where
\[
\bar{J}_n(x):=\frac{1-e^{-n(\lossexpn'(b)-x)\delta}}{n(\lossexpn'(b)-x)}e^{-nd_n(b)}.
\]
At this stage, denote
\[
{J}_n(x):=\frac{e^{-nd_n(b)}}{n(\lossexpn'(b)-x)}.
\]
Fix $x<(\te_0-b)f'(b)/\{2(1-b)\}$, so that $x < \lossexpn'(b)$ for sufficiently large $n$.
If $n$ is sufficiently large, then $(1-\vep)J_n(x)<\bar{J}_n(x)<J_n(x)$.
Hence, \eqref{f: let's integrate !!} becomes
\begin{equation}\label{f: let's integrate !!!}
(p(b)-\vep)\,(1-\vep)\,J_n(-\vep)
\leq J_n \leq
(p(b)+\vep)\, J_n(\vep)\,+\,
k\,e^{-n\lossexpn(b+\delta)}.
\end{equation}

Since $\lossexpn$ is increasing on $(b+\delta,\:1)\subset (\te_0,\:1)$,
$e^{-n\lossexpn(b+\delta)}=o(J_n(x))$ as $n\to \infty$ for
$x<(\te_0-b)f'(b)/\{2(1-b)\}$.
In fact, $\lossexpn(b)$, $\lossexpn'(b)$ and $\lossexpn(b+\delta)$ converge to positive constants
by the strong law of large numbers,
$f$ (and therefore $\lossexpn$ and $\lossexpn'$) being continuous.
Hence,
\eqref{f: let's integrate !!!} yields
\[(p(b)-\vep)(1-\vep)\,J_n(-\vep)
\leq J_n \leq
(p(b)+2\vep)\, J_n(\vep)\]
for sufficiently large $n$. The number $\vep$ being arbitrary, it follows that
\[J_n\sim p(b)\,J_n(0)\]
as $n\to \infty$ and \eqref{f: Num2} is proved.

In order to prove \eqref{f: Num1}, take
$\bar{\lossexpn}(t)\,:=\,d(t,\,1-\mle)\,=\,\lossexpn(1-t)$,
$\bar{p}(t):=p(1-t)$, $\bar{\te}_0:=1-\te_0$
$\bar{b}:=1-a$ (so that $\bar{b}>\bar{\te}_0$) and 
notice that by \eqref{f: Num2}
\begin{equation*}
\int_{(\bar{b},\: 1)} \bar{p}(s)\, e^{-n\,\bar{\lossexpn}(s)}\, \diff s \sim
\frac{1}{n}\, \frac{\bar{p}(\bar{b})}{\bar{\lossexpn}'(\bar{b})}\,
e^{-n\,\bar{\lossexpn}(\bar{b})},
\end{equation*}
and then apply the substitution $t=1-s$ in the integral.

\end{proof}

\begin{proof}[Proof of Proposition \ref{prop: conv.rate}]
Let $p=\diff \pi/\diff \la$, where $\la$ is the Lebesgue measure.
%Since \eqref{f: support} is true for every possible value for $\te_0$, the
%density $p$ must be positive almost everywhere. Hence, Lemma \ref{lemma: var.norm} can be applied.
%Notice that
In order to apply Lemma \ref{lemma: var.norm}, notice that
\[ \post((\te_0-\vep,\te_0+\vep)^c) =
\frac{\int_{(0,\te_0-\vep)}
e^{-n\, d(\te,\mle)} p(\te)\,\diff \te+
\int_{(\te_0+\vep,1)}
e^{-n\, d(\te,\mle)} p(\te)\,\diff \te}
{\int_{(0,1)}
e^{-n\, d(\te,\mle)} p(\te)\,\diff \te}. \]
Hence, combining \eqref{f: Den} with \eqref{f: Num1} and \eqref{f: Num2}, the thesis follows.

\end{proof}

\begin{proof}[Proof of Proposition \ref{prop: variance}]

Denote by $p(\cdot)$ the density of $\pi$ with respect to the Lebesgue measure.
To begin, notice that
\begin{equation}\label{f: var}
\begin{split}
\int_{(0,1)}(t-\mle)^2\, \post(\diff t)&=
\frac{\int_{(0,1)}\, e^{-n \, d(t,\mle)}\,(t-\mle)^2\, \diff \pi(t)}
{\int_{(0,1)} e^{-n \, d(t,\mle)}\, \diff \pi(t)} \\
&=\frac{\int_{(0,1)}p(t)\, e^{-n \, d(t,\mle)}\,(t-\mle)^2\, \diff t}
{\int_{(0,1)} p(t)\,e^{-n \, d(t,\mle)}\, \diff t} .
\end{split}\end{equation}
%Since \eqref{f: support} is true for every possible value for $\te_0$, the density $p$ must be positive almost everywhere.
Lemma \ref{lemma: var.norm} can be applied for both
the numerator and the denominator of \eqref{f: var}.
Combination of \eqref{f: var} with \eqref{f: Den} and \eqref{f: VAR}
yields
\begin{equation}\label{f: approx.norm.}
\int_{(0,1)}(t-\mle)^2\, \post(\diff t)
\sim \frac{1}{n\,\lossexpn''(\mle)},%\sim \frac{1}{n\,\lossexpn''(\te_0)} ,
\end{equation}
as $n\to \infty$.

Combining \eqref{f: approx.norm.} with \eqref{f: second.deriv..},
one obtains that
\begin{equation}\label{f: cond.var}
\int_{(0,1)}(t-\mle)^2\, \post(\diff t)\sim
-\frac{1-\mle}{n\,f'(\mle)},
\end{equation}
as $n\to \infty$.

In virtue of continuity of $f'$, by the strong law of large numbers,
\eqref{f: cond.var} entails that
\begin{equation}\label{f: cond.var..}
\int_{(0,1)}(t-\mle)^2\,
\post(\diff t)\sim -\frac{1-\te_0}{nf'(\te_0)},
\end{equation}
as $n\to \infty$. %, $P_{\te_0}$ - a.s.

%Let $T_n:=T_n(\mle,f)$ be a random variable valued into $(0,\, 1)$ with distribution $\post$.
Let $E_f(\mle)$ denote the the mean with respect to the distribution $\post$. 

Notice that
\begin{equation}\label{f: variances1}
V_f(\mle)=\int_{(0,1)}(t-\mle)^2\, \post(\diff t)\,-\,(E_f(\mle)-\mle)^2,
%\var(T_n)= \mean(T_n-\mle)^2\,-\,(\mean(T_n)-\mle)^2
\end{equation}
and
\begin{equation}\label{f: E_f}
\int_{(0,1)}(t-\mle)^2\, \post(\diff t)
= \int_{(0,1)}t^2\, \post(\diff t)\,-\,2\mle\int_{(0,1)}t\, \post(\diff t)\,+\,\mle^2.
\end{equation}
By \eqref{f: E_f}, we obtain that
\begin{equation}\label{f: XXXY}\begin{split}
\mod{E_f(\mle)-\mle}&=\smallfrac{1}{2\mle}
\mod{\int_{(0,1)}t^2\, \post(\diff t) -\mle^2 - \int_{(0,1)}(t-\mle)^2\, \post(\diff t)}\\
&\leq
\smallfrac{1}{2\mle}
\bigg(\mod{\int_{(0,1)}t^2\, \post(\diff t) -\mle^2} +
\int_{(0,1)}(t-\mle)^2\, \post(\diff t)\bigg)
\end{split}\end{equation}
At this stage, dividing \eqref{f: VARR} by \eqref{f: Den}
and applying \eqref{f: second.deriv..} and
the strong law of large numbers, one obtains that
\begin{equation}\label{f: AAA}
\int_{(0,1)}t^2\, \post(\diff t)- \mle^2\sim -\frac{1-\te_0}{nf'(\te_0)}.
\end{equation}
In virtue of \eqref{f: cond.var..} and \eqref{f: AAA},
equation \eqref{f: XXXY} yields:
\begin{equation*}%\label{f: second.m...}
E_f(\mle)-\mle=O\Big(\frac{1}{n}\Big).
\end{equation*}
Hence, $(E_f(\mle)-\mle)^2$ is negligible with respect to \eqref{f: cond.var..} and
therefore \eqref{f: variances1} entails that
\begin{equation}\label{f: almost.there}
V_f(\mle)\sim\int_{(0,1)}(t-\mle)^2\, \post(\diff t)%\qquad P_{\te_0}-\text{a.s.}
\end{equation}
The thesis follows from \eqref{f: almost.there} and \eqref{f: cond.var..}.
\end{proof}

\section*{Acknowledgements}
We thank two anonymous referees for their helpful comments. 
This work was partially supported by ESF and Regione Lombardia (by the grant ``Dote Ricercatori'').

%%%\bibliographystyle{elsarticle-harv}
%\bibliographystyle{agsm}
%\bibliography{Bibliogr}

\end{document}